\documentclass{tran-l}
\usepackage{amssymb} 
\usepackage{eucal}   

\renewcommand\ge\geqslant
\renewcommand\geq\geqslant
\renewcommand\le\leqslant
\renewcommand\leq\leqslant

\theoremstyle{plain}
\newtheorem{theorem}{Theorem}[section]

\theoremstyle{definition}
\newtheorem*{definition}{Definition}

\theoremstyle{remark}
\newtheorem*{remark}{Remark}

\numberwithin{equation}{section}

\makeatletter
\newenvironment{sizeequation}[1]{%
  \skip@=\baselineskip
  #1%
  \baselineskip=\skip@
  \equation
}{\endequation \ignorespacesafterend}
\newenvironment{sizedisplaymath}[1]{%
  \skip@=\baselineskip
  #1%
  \baselineskip=\skip@
  \displaymath
}{\enddisplaymath \ignorespacesafterend}
\newenvironment{sizeeqnarray}[1]{%
  \skip@=\baselineskip
  #1%
  \baselineskip=\skip@
  \eqnarray
}{\endeqnarray \ignorespacesafterend}

\newenvironment{sizeeqnarray*}[1]{%
  \def\@eqncr{\nonumber\@seqncr}
  \skip@=\baselineskip
  #1%
  \baselineskip=\skip@
  \eqnarray
}{\nonumber\endeqnarray \ignorespacesafterend}
\makeatother

\DeclareMathOperator{\angliceRe}{Re}
\newcommand{\tusp}{{$\,$}}

\begin{document}

\title[A Generalization of Euler's Hypergeometric Transformation]{A Generalization of Euler's\\ Hypergeometric Transformation}
\author{Robert S. Maier}
\address{Departments of Mathematics and Physics,
University of Arizona, Tucson, AZ 85721, USA}
\thanks{Partially supported by NSF grant PHY-0099484.}
\email{rsm@math.arizona.edu}
\urladdr{http://www.math.arizona.edu/{\textasciitilde}rsm}
\subjclass{Primary 33C20; Secondary 33C05, 34Mxx}
\copyrightinfo{2005}{American Mathematical Society}

\begin{abstract}
Euler's transformation formula for the Gauss hypergeometric
function~${}_2F_1$ is extended to hypergeometric functions of higher order.
Unusually, the generalized transformation constrains the hypergeometric
function parameters algebraically but not linearly.  Its consequences for
hypergeometric summation are explored.  It~has as corollary a summation
formula of Slater.  From this formula new one-term evaluations of
${}_2F_1(-1)$ and ${}_3F_2(1)$ are derived, by applying transformations in
the Thomae group.  Their parameters are also constrained nonlinearly.
Several new one-term evaluations of ${}_2F_1(-1)$ with linearly constrained
parameters are derived as~well.
\end{abstract}

\maketitle

\section{Introduction}

Many hypergeometric identities are known (see Slater~\cite{Slater66} and
Prudnikov et~al.\ \cite[Chapter~7]{Prudnikov86c}).  Most are summation
formulas for hypergeometric series, either infinite or terminating,
including binomial coefficient summations \cite{Gould72}.  Recent work has
concentrated on developing techniques for verifying asserted or conjectured
identities, rather on than deriving new ones.  Zeilberger's algorithm and
the Wilf--Zeilberger method can verify nearly all classical hypergeometric
identities, at~least those involving terminating series, and many more
besides~\cite{Koepf98,Petkovsek96}.  Under the circumstances, maintaining
an exhaustive list of such identities has come to be seen as unnecessary.
It~has even been provocatively asserted that ``there is no hypergeometric
database,'' in the sense that the class of hypergeometric identities is for
practical purposes open-ended, and that such identities are best proved by
algorithmic techniques more sophisticated than database lookup
\cite[Chapter~3]{Petkovsek96}.

But hypergeometric {\em function transformations\/}, in which the argument
is a free variable, are more manageable than general identities.  One could
reasonably hope to enumerate or otherwise characterize the class, say, of
all two-term function transformations.  In~fact the two-term
transformations relating the Gauss hypergeometric function ${}_2F_1(x)$ to
${}_2F_1(Rx)$, where $R$~is a rational map of the Riemann sphere to itself,
are now fully classified.  Besides the celebrated transformation of Euler,
\begin{align}
{}_2F_1\biggl(
\begin{array}{c}
{a,\,b} \\
{c}\\
\end{array}
\!
\biggm|
x
\biggr)
&=
(1-x)^{c-a-b}\ 
{}_2F_1\biggl(
\begin{array}{c}
{c-a,\,c-b} \\
{c}\\
\end{array}
\!
\biggm|
x
\biggr),
\label{eq:euler}
\end{align}
in which $R(x)=x$, and Pfaff's transformation, in which $R(x)=x/(x-1)$,
there are transformations of larger mapping degree ($\deg R>1$), which were
classified by Goursat~\cite{Goursat1881}.  The best known are the quadratic
ones, which were originally worked~out by Gauss and Kummer, and proved
concisely by Riemann~\cite[\S\tusp3.9]{Andrews99}.  Recently, Goursat's
classification was completed by an enumeration of the transformations
of~${}_2F_1$ with no free parameter, most of which have quite large
degree~\cite{Vidunas2004}.  Several of the quadratic and cubic
transformations of~${}_2F_1$ have analogues on the ${}_3F_2$ level,
discovered by Whipple and Bailey respectively~\cite{Askey94}.  `Companion'
transformations relate ${}_3F_2$ to ${}_4F_3$~\cite{Gessel82}.  But
no~clear analogues on levels above~${}_2F_1$ of the remaining
${}_2F_1$~transformations, in particular of the degree-$1$ transformations
of Euler and Pfaff, have previously been found.

We report here on a recent discovery: Euler's
transformation~(\ref{eq:euler}) has an analogue on all higher levels.  That
this has not been noticed before may be due~to the fact that in the
generalized transformation, the hypergeometric parameters are constrained
{\em algebraically\/}, not linearly.  This is very unusual.  In~the
compilation of Prudnikov et~al.\ \cite[Chapter~7]{Prudnikov86c}, no
hypergeometric identity with parameters constrained in this way is listed.
The only published hypergeometric function transformation with nonlinearly
constrained parameters that we have been able to uncover is a quadratic one
connecting ${}_3F_2$ and ${}_5F_4$ \cite{Niblett51}, which seems unrelated.

The generalized Euler's transformation came close to being discovered in
the 1950s by Slater, who obtained a formula for ${}_{r+1}F_r(1)$ that
applies when its parameter vector is restricted to a certain algebraic
variety~\cite{Slater55}.  Several closely related summation formulas with
nonlinear parametric restrictions had been obtained much earlier by Searle
\cite{Searle09}.  Searle's and Slater's formulas have not attracted the
attention they deserve, though Slater's was later reproduced in
\cite[\S\tusp2.6.1]{Slater66}.  More recently, a few additional expressions
for hypergeometric sums with algebraically constrained parameters have been
published.  These include a formula for a terminating ${}_4F_3(1)$
\cite[\S\tusp3.13]{Luke69}, a related one for a terminating ${}_3F_2(1)$
\cite[Eq.~(1.9)]{Gessel82}, an evaluation of a non-terminating ${}_4F_3(1)$
\cite{Gosper76}, and an exotic evaluation of a non-terminating
${}_2F_1(-1)$ \cite[\S\tusp4]{Vidunas2001a}.

This article is organized as~follows.  The generalized Euler's
transformation appears as Theorem~\ref{thm:main}.  A~combinatorial proof is
given, and Slater's summation formula and several additional hypergeometric
identities are derived as corollaries.  In Section~\ref{sec:3F2} the
implications for ${}_3F_2(1)$ are examined.  Using its $S_5$~symmetry, a
set of three related ${}_3F_2(1)$ evaluations is derived, each with
nonlinear parametric constraints and three free parameters.  The extent to
which they overlap the standard identities of Dixon, Watson, Whipple, and
Pfaff--Saalsch\"utz is determined.  In Section~\ref{sec:2F1} these new
identities are employed to generate new evaluations of ${}_2F_1(-1)$,
including several exotic ones with nonlinearly constrained parameters.

\section{Key results}
\label{sec:key}

The following standard notation will be used.  For any integer $r\ge1$ and
parameter vector
$(\alpha;\beta)=(\alpha_1,\dots,\alpha_{r+1};\beta_1,\dots,\beta_r)\in\mathbb{C}^{r+1}\times\mathbb{C}^r$
in which no~$\beta_i$ is a non-positive integer, the generalized
hypergeometric series is
\begin{equation}
\label{eq:hyperseries}
{}_{r+1}F_r
\biggl(
\biggl.
\begin{array}{c}
{\alpha_1,\,\ldots,\,\alpha_{r+1}} \\
{\beta_1,\,\ldots,\,\beta_r}\\
\end{array}
\biggr| 
\ x
\biggr)
=
\sum_{k=0}^\infty
\frac
{(\alpha_1)_k\cdots(\alpha_{r+1})_k}
{(\beta_1)_k\cdots(\beta_{r})_k}
\ \frac{x^k}{k!},
\end{equation}
where $(\alpha)_k$ signifies the rising factorial
$\alpha(\alpha+1)\cdots(\alpha+k-1)$, with ${(\alpha)_0=1}$.  If $x$ is
omitted, $x=1$ is understood.  The series converges absolutely on~$|x|<1$,
and if the {\em parametric excess\/} $s=\sum_{i=1}^r\beta_i -
\sum_{i=1}^{r+1}\alpha_i$ has positive real part, it converges at~$x=1$.
(For convergence at~$x=-1$, $\angliceRe s>-1$ suffices.)  The series can be
continued uniquely from the unit disk to the cut Riemann sphere
$\mathbb{CP}^1\setminus[1,\infty]$.  Any hypergeometric transformation
formula, when stated without restriction on~$x$, should be taken as holding
in a neighborhood of~$x=0$, namely the largest neighborhood of~$x=0$
in~$\mathbb{CP}^1\setminus[1,\infty]$ to which both sides can be continued.

A generalized hypergeometric series or function is said to be {\em
$s$-balanced\/} if the parametric excess equals~$s$.  It is {\em
well-poised\/} if the parameters $\{\alpha_i\},\{\beta_i\}$ can be
separately permuted so that $\alpha_1 + 1 = \alpha_2 + \beta_1 = \cdots =
\alpha_{r+1} + \beta_r$.

The generalized Euler transformation formula will now be stated and proved.

\begin{definition}
For each $r\ge1$, the algebraic variety
${\mathcal{U}}_r\subset\mathbb{C}^{r+1}\times\mathbb{C}^r$, which is
$(r+1)$-dimensional, comprises all~$(\alpha;\beta)$ for which the
$r$~equations
\begin{displaymath}
\left\{
\begin{array}{rcl}
\sum_{1\le i\le r+1}\alpha_i&=&\sum_{1\le i\le r}\beta_i  \\
\sum_{1\le i<j\le r+1}\alpha_i\alpha_j&=&\sum_{1\le i<j\le r}\beta_i\beta_j \\
\sum_{1\le i<j<k\le r+1}\alpha_i\alpha_j\alpha_k&=&\sum_{1\le i<j<k\le r}\beta_i\beta_j\beta_k \\
&\vdots&
\end{array}
\right.
\end{displaymath}
are satisfied.  In the $k$'th equation, the left and right sides are the
$k$'th elementary symmetric polynomial in $\alpha_1,\dots,\alpha_{r+1}$ and
in~$\beta_1,\dots,\beta_{r}$, respectively.  
\end{definition}

\begin{theorem}
For all $r\ge1$ and $(\alpha;\beta)\in\mathcal{U}_r$ for which
no~$\beta_i+1$ is a non-positive integer,
{\small
\begin{displaymath}
{}_{r+1}F_r\biggl(
\begin{array}{c}
{\alpha_1,\,\ldots,\,\alpha_{r+1}} \\
{\beta_1+1,\,\ldots,\,\beta_r+1}\\
\end{array}
\!
\biggm|
x
\biggr)
=
(1-x)\ 
{}_{r+1}F_r\biggl(
\begin{array}{c}
{\alpha_1+1,\,\ldots,\,\alpha_{r+1}+1} \\
{\beta_1+1,\,\ldots,\,\beta_r+1}\\
\end{array}
\!
\biggm|
x
\biggr).
\end{displaymath}
}
\label{thm:main}
\end{theorem}

\begin{remark}
By the condition $(\alpha;\beta)\in\mathcal{U}_r$, the series on the left
and right sides are respectively $r$-balanced and $(-1)$-balanced.  If
$r=1$, this reduces to the $1$-balanced case of Euler's transformation,
i.e., the case $c-a-b=1$ of~(\ref{eq:euler}).
\end{remark}

\begin{proof}[Proof of Theorem~\ref{thm:main}]
By the series representation~(\ref{eq:hyperseries}), the theorem will
follow if for all $k\ge1$,
\begin{sizeeqnarray*}
{\small}
&&\frac
{(\alpha_1)_k\,\cdots\,(\alpha_{r+1})_k}
{(\beta_1+1)_k\,\cdots\,(\beta_{r}+1)_k \,k!}
\\
&&\qquad\qquad
=\frac
{(\alpha_1+1)_k\,\cdots\,(\alpha_{r+1}+1)_k}
{(\beta_1+1)_k\,\cdots\,(\beta_{r}+1)_k \,k!}
\,-\,
\frac
{(\alpha_1+1)_{k-1}\,\cdots\,(\alpha_{r+1}+1)_{k-1}}
{(\beta_1+1)_{k-1}\,\cdots\,(\beta_{r}+1)_{k-1}\,(k-1)!}.
\end{sizeeqnarray*}
By examination, this equation can be obtained from
\begin{equation}
\label{eq:needed2}
\alpha_1\cdots\alpha_{r+1}
=
(\alpha_1+k)\cdots(\alpha_{r+1}+k)
-
(\beta_1+k)\cdots(\beta_{r}+k)k
\end{equation}
by multiplying by the product
$(\alpha_1+1)_{k-1}\cdots(\alpha_{r+1}+1)_{k-1}$ and dividing by the
product $(\beta_1+1)_k\cdots(\beta_r+1)_k k!$.  The left side
of~(\ref{eq:needed2}) is independent of~$k$, and the right side is a
polynomial in~$k$ of degree~$r$.  For $n=1,\dots,r$, the coefficient
of~$k^n$ is proportional to the sum of all monomials
$\alpha_{i_1}\cdots\alpha_{i_{r+1-n}}$ minus the sum of all monomials
$\beta_{j_1}\cdots\beta_{j_{r+1-n}}$.  But since
$(\alpha;\beta)\in\mathcal{U}_r$, each coefficient is zero.  So the right
side of~(\ref{eq:needed2}) reduces to a constant, namely
$\alpha_1\cdots\alpha_{r+1}$, and is the same as the left side.  The
theorem is proved.
\end{proof}

Theorem~\ref{thm:main} can be viewed as belonging to the theory of
contiguous function relations, which (for~${}_2F_1$) dates back to Gauss.
Recall that if a parameter of~${}_{r+1}F_r$ is displaced by~$\pm1$, the
resulting function is said to be contiguous to the original.  Gauss showed
that ${}_2F_1(a,b;c;x)$ and any two of its contiguous series are connected
by a homogeneous linear relation with coefficients linear in~$x$ and
polynomial in~$(a,b;c)$.  In~fact by constructing chains of contiguous
functions, he showed that any series of the form ${}_2F_1(a+l,b+m;c+n;x)$,
where $(l,m;n)\in\mathbb{Z}^2\times\mathbb{Z}$, is connected to
${}_2F_1(a,b;c;x)$ and any one of its contiguous series by a three-term
homogeneous linear relation, the coefficients of which are polynomial in
$x$ and $(a,b;c)$.  By imposing parametric constraints, i.e., confining
$(a,b;c)$ to an algebraic variety, it is sometimes possible to obtain
two-term relations, as~well.

The two-term relation of Theorem~\ref{thm:main} is clearly related to this
classical theory, but is much stronger in~that it holds for all~$r\ge1$.
Contiguous function relations for arbitrary ${}_{r+1}F_r$ (and indeed for
arbitrary~${}_pF_q$) were worked~out by Rainville~\cite{Rainville45}.
There is a fundamental set of $3r+2$ such relations, with up~to $r+2$
terms, from which all others follow by linearity.  But the theory of those
relations is complicated, since Gauss's construction does not extend
to~$r>1$.  In particular, no two-term relation for general~$r$ was
previously known.  Theorem~\ref{thm:main} is therefore surprising.

When applied to Euler's transformation, the technique of equating
coefficients produces the Pfaff--Saalsch\"utz formula, which sums any
terminating $1$-balanced ${}_3F_2(1)$~\cite[\S\tusp2.2]{Andrews99}.
Applied to Theorem~\ref{thm:main}, it produces the following.

\begin{theorem}
\label{thm:terminatingslater}
For all $r\geq1$ and $(\alpha;\beta)\in\mathcal{U}_r$ for which
no $\beta_i+1$ is a non-positive integer, 
\begin{equation}
\label{eq:newnew}
\sum_{j=0}^n
\frac
{(\alpha_1)_j\,\cdots\,(\alpha_{r+1})_j}
{(\beta_1+1)_j\,\cdots\,(\beta_{r}+1)_j\,j!}
=
\frac
{(\alpha_1+1)_n\,\cdots\,(\alpha_{r+1}+1)_n}
{(\beta_1+1)_n\,\cdots\,(\beta_{r}+1)_n\,n!}
\end{equation}
for every integer $n\ge0$.  In consequence, for all $r\geq1$ and
$(\alpha;\beta)\in\mathcal{U}_r$ for which $\alpha_1$ equals~$-n$, a
non-positive integer, and no $\beta_i+1$ is a non-positive integer,
\begin{equation}
\label{eq:newclass}
{}_{r+1}F_r 
\biggl(
\begin{array}{c}
{-n,\,\alpha_2,\,\ldots,\,\alpha_{r+1}} \\
{\beta_1+1,\,\ldots,\,\beta_r+1}\\
\end{array}
\biggr)
=
\frac
{(\alpha_1+1)_n\,\cdots\,(\alpha_{r+1}+1)_n}
{(\beta_1+1)_n\,\cdots\,(\beta_{r}+1)_n\,n!}
\end{equation}
\end{theorem}

\begin{proof}
To obtain~(\ref{eq:newnew}) from Theorem~\ref{thm:main}, multiply both
sides by $(1-x)^{-1}$ and equate the coefficients of~$x^n$ on the two
sides.
\end{proof}

The formula~(\ref{eq:newclass}) sums a class of terminating $r$-balanced
${}_{r+1}F_r(1)$ series, in which the parameters are required to satisfy
$r-1$ nonlinear conditions.  It~is worth considering the relation between
(\ref{eq:newclass}) and other summation formulas.  Any terminating
$2$-balanced ${}_3F_2(1)$ can be summed by the Sheppard--Anderson
formula~\cite{Roy87}, so the $r=2$ case of~(\ref{eq:newclass}) is a special
case of a more general known result.  When $r>2$, the situation is quite
different, since few closed-form summations of terminating $r$-balanced
${}_{r+1}F_r(1)$'s with $r>2$ are known.  Certain terminating $4$-balanced
${}_5F_4(1)$'s~\cite[\S\tusp6]{Stanton98} and terminating $6$-balanced
${}_7F_6(1)$'s~\cite{Askey89} can be evaluated in closed form, though
explicit formulas have not been published.  It~would be interesting to
determine the extent to which (\ref{eq:newclass}) overlaps those
closed-form evaluations.  For larger values of~$r$, (\ref{eq:newclass})~is
in a league by itself.

The following theorem is the extension of~(\ref{eq:newclass}) to
non-terminating series.

\begin{theorem}
\label{thm:slater}
For all $r\geq1$ and $(\alpha;\beta)\in\mathcal{U}_r$ for which no
$\beta_i+1$ is a non-positive integer,
\begin{displaymath}
{}_{r+1}F_r 
\biggl(
\begin{array}{c}
{\alpha_1,\,\ldots,\,\alpha_{r+1}} \\
{\beta_1+1,\,\ldots,\,\beta_r+1}\\
\end{array}
\biggr)
=
\Gamma
\biggl[
\begin{array}{c}
{\beta_1+1,\,\ldots,\,\beta_{r}+1} \\
{\alpha_1+1,\,\ldots,\,\alpha_{r+1}+1}\\
\end{array}
\biggr].
\end{displaymath}
\end{theorem}

\begin{remark}
This introduces a useful notation.  The right side is a gamma quotient, of
the sort that appears in many hypergeometric summations.  It~would be
written in~full as $\Gamma(\beta_1+1)\cdots\Gamma(\beta_r+1)/
\Gamma(\alpha_1+1)\cdots\Gamma(\alpha_{r+1}+1)$.  Any gamma quotient with a
denominator argument equal to a non-positive integer is interpreted as
zero.
\end{remark}

\begin{proof}[Proof of Theorem~\ref{thm:slater}]
Take $n\to\infty$ in~(\ref{eq:newnew}), using the $p,q=r+1,r$ case of
\begin{sizedisplaymath}
{\small}
\frac{(\alpha_1+1)_n\,\cdots\,(\alpha_{p}+1)_n}
{(\beta_1+1)_n\,\cdots\,(\beta_{q}+1)_n\,n!}
\sim \frac{\Gamma(\beta_1+1)\,\cdots\,\Gamma(\beta_{q}+1)}
{\Gamma(\alpha_1+1)\,\cdots\,\Gamma(\alpha_{p}+1)}\ 
n^{\sum_{i=1}^p(\alpha_i+1)-\sum_{i=1}^q(\beta_i+1) - 1},
\end{sizedisplaymath}
which is well-known $n\to\infty$ asymptotic
behavior~\cite[\S\tusp2.1]{Andrews99}.
\end{proof}

Theorems~\ref{thm:terminatingslater} and~\ref{thm:slater} are precisely the
summation results of Slater.  The best reference for her original
derivation is the note~\cite{Slater55}, since the version in
\cite[\S\tusp2.6.1]{Slater66} contains several unfortunate transcription
errors.  In~particular, the quantity $a+b+c-d$ in (2.6.1.10)--(2.6.1.11)
should be read as $2+a+b+c-d$.  It~should be mentioned that much earlier,
Searle~\cite{Searle09} proved a related theorem on partial sums of series
of hypergeometric type, from which Slater's results can be deduced as
corollaries.  (The author is indebted to George Andrews for this
reference.)

For every integer $m\ge0$, there is a similar summation formula for a class
of $(r+m)$-balanced ${}_{r+1}F_r(1)$ series, with the parameters required
to satisfy $r-1$ nonlinear conditions.  This result is best framed in~terms
of an $(r+1)$-dimensional algebraic variety
$\mathcal{U}_r^m\subset\mathbb{C}^{r+1}\times\mathbb{C}^r$, indexed by~$m$,
with $\mathcal{U}_r^0$ equalling $\mathcal{U}_r$.  The $r$~polynomial
equations defining $\mathcal{U}_r^m$ become unwieldy when $m$~is large, so
the following exposition of the $m=1$ case should suffice.

If $(\alpha;\beta)\in\mathbb{C}^{r+1}\times\mathbb{C}^r$, let
$S_l(\alpha)$, resp.\ $S_l(\beta)$, denote the sum of the $\binom{r+1}{l}$
monomials of the form $\alpha_{i_1}\cdots\alpha_{i_l}$, resp.\ the
$\binom{r}{l}$ monomials of the form $\beta_{j_1}\cdots\beta_{j_l}$, with
$S_0(\alpha)=S_0(\beta)=1$ by convention.  $\mathcal{U}_r^0=\mathcal{U}_r$
is the common solution set of the equations $S_l(\alpha)-S_l(\beta)=0$,
$l=1,\dotsc,r$.  Let
$\mathcal{U}_r^1\subset\mathbb{C}^{r+1}\times\mathbb{C}^r$ be the common
solution set of the equations
\begin{eqnarray}
\label{eq:newslatervariety}
&&\left[S_{l-1}(\alpha)-S_{l-1}(\beta)\right]
\left[S_{r+1}(\alpha)+S_{r}(\beta)\right]\\
&&\qquad{}+
\left[S_{l}(\alpha)-S_{l}(\beta)+S_{l-1}(\beta)\right]
\left[-S_{r}(\alpha)+S_{r}(\beta)\right]
=0,
\nonumber
\end{eqnarray}
$l=1,\dotsc,r$.  It follows from the $l=1$ equation that
$(\alpha;\beta)\in\mathcal{U}^1_r$ only if the series
${}_{r+1}F_r(\alpha_1,\dotsc,\alpha_{r+1};\beta_1+1,\dotsc,\beta_r+1;1)$ is
$(r+1)$-balanced, unless $S_r(\alpha)=S_r(\beta)$.

\begin{theorem}
\label{thm:slaterplus}
For all $r\geq1$ and $(\alpha;\beta)\in\mathcal{U}_r^1$ for which no
$\beta_i+1$ is a non-positive integer and for which $S_{r+1}(\alpha)+
S_r(\beta)\neq0$,
\begin{displaymath}
{}_{r+1}F_r 
\biggl(
\begin{array}{c}
{\alpha_1,\,\ldots,\,\alpha_{r+1}} \\
{\beta_1+1,\,\ldots,\,\beta_r+1}\\
\end{array}
\biggr)
=
\left[
\frac
{-S_{r}(\alpha)+S_r(\beta)}
{S_{r+1}(\alpha)+S_r(\beta)}
\right]
\,
\Gamma
\biggl[
\begin{array}{c}
{\beta_1+1,\,\ldots,\,\beta_{r}+1} \\
{\alpha_1+1,\,\ldots,\,\alpha_{r+1}+1}\\
\end{array}
\biggr].
\end{displaymath}
\end{theorem}

\begin{proof}
Trivially,
\begin{displaymath}
{}_{r+1}F_r 
\biggl(
\begin{array}{c}
{\alpha_1,\,\ldots,\,\alpha_{r+1}} \\
{\beta_1+1,\,\ldots,\,\beta_r+1}\\
\end{array}
\biggr)
=
{}_{r+2}F_{r+1}
\biggl(
\begin{array}{c}
{\alpha_1,\,\ldots,\,\alpha_{r+1},\,\gamma} \\
{\beta_1+1,\,\ldots,\,\beta_r+1,\,\gamma}\\
\end{array}
\biggr).
\end{displaymath}
If the conditions $S_l(\alpha_1,\dotsc,\alpha_{r+1},\gamma) -
S_l(\beta_1,\dotsc,\beta_{r},\gamma)=0$, $l=1,\dotsc,r+1$, are satisfied,
then the right-hand series can be summed by Theorem~\ref{thm:slater}.  It
follows from the condition labeled by $l=r+1$ that for this to occur,
$\gamma$~must equal
$\left[S_{r+1}(\alpha)+S_r(\beta)\right]/\left[-S_r(\alpha)+S_r(\beta)\right]$.
Eliminating~$\gamma$ from the remaining conditions
yields~(\ref{eq:newslatervariety}).  The prefactor in the theorem is simply
$\Gamma(\gamma)/\Gamma(\gamma+1)=1/\gamma$.
\end{proof}

So a class of $(r+1)$-balanced ${}_{r+1}F_r(1)$'s with nonlinearly
constrained parameters can be summed in closed form.  This is the $m=1$
case of the result mentioned above.  The extension to $m>1$ is
algebraically nontrivial.  For $m=2$, one can write
\begin{displaymath}
{}_{r+1}F_r 
\biggl(
\begin{array}{c}
{\alpha_1,\,\ldots,\,\alpha_{r+1}} \\
{\beta_1+1,\,\ldots,\,\beta_r+1}\\
\end{array}
\biggr)
=
{}_{r+3}F_{r+2}
\biggl(
\begin{array}{c}
{\alpha_1,\,\ldots,\,\alpha_{r+1},\,\gamma,\,\delta} \\
{\beta_1+1,\,\ldots,\,\beta_r+1,\,\gamma,\,\delta}\\
\end{array}
\biggr),
\end{displaymath}
and solve for the $\gamma,\delta$ for which the right-hand side can be
evaluated by Theorem~\ref{thm:slater}.  
Each of $\gamma,\delta$ will
be an irrational algebraic function of~$(\alpha;\beta)$, and the same will
be true of the prefactor multiplying the gamma quotient.  In~fact, this
prefactor will be an irrational algebraic function of~$(\alpha;\beta)$ for
all~$m>1$.  So the statement of the extension of
Theorem~\ref{thm:slaterplus} to $m>1$, $(\alpha;\beta)\in\mathcal{U}_r^m$,
is rather complicated.

Theorems \ref{thm:terminatingslater}--\ref{thm:slaterplus} followed from
Theorem~\ref{thm:main}, the generalized Euler's transformation, by the
method of equating coefficients.  The beta integral method, which is a more
sophisticated way of generating hypergeometric summation formulas from
transformation formulas, can also be used.  Krattenthaler and Srinivasa
Rao~\cite{Krattenthaler2003} have applied this method to the known
transformations of ${}_{r+1}F_r(x)$ to ${}_{r+1}F_r(Rx)$.  Applying it to
Theorem~\ref{thm:main} yields the following remarkable identity.

\begin{theorem}
\label{thm:genThomae}
For all $r\ge1$ and $(\alpha;\beta)\in\mathcal{U}_r$ for which
no~$\beta_i+1$ is a non-positive integer,
\begin{multline*}
{}_{r+2}F_{r+1}\biggl(
\begin{array}{c}
{\alpha_1,\,\ldots,\,\alpha_{r+1},\,A} \\
{\beta_1+1,\,\ldots,\,\beta_r+1,\,B}\\
\end{array}
\biggr)\\
=
\left(\frac{B-A}B\right)
{}_{r+2}F_{r+1}\biggl(
\begin{array}{c}
{\alpha_1+1,\,\ldots,\,\alpha_{r+1}+1,\,A} \\
{\beta_1+1,\,\ldots,\,\beta_r+1,\,B+1}\\
\end{array}
\biggr),
\end{multline*}
for all $A,B\in\mathbb{C}$, provided\/ $\angliceRe(B-A)>0$ and $B$~is
not a non-positive integer.
\label{thm:appendage}
\end{theorem}

\begin{remark}
Since $(\alpha;\beta)\in\mathcal{U}_r$, the two series are respectively
$(B-A+r)$-balanced and $(B-A)$-balanced.  If $r=1$ there are no nonlinear
constraints, and this identity reduces to a special case of one of Thomae's
transformations of ${}_3F_2(1)$ (see below).
\end{remark}

\begin{proof}
[Proof of Theorem~\ref{thm:appendage}] Multiply both sides of
Theorem~\ref{thm:main} by $x^{A-1}(1-x)^{B-A-1}$, and integrate from $x=0$
to~$x=1$.  To evaluate each side, interchange integration and summation,
and use the fact that the beta integral $\smallint_0^1
x^{\mu-1}(1-x)^{\nu-1}\,dx$ equals
$\Gamma(\mu)\Gamma(\nu)/\Gamma(\mu+\nu)$.  Finally, convert each side back
to hypergeometric notation.
\end{proof}


\section{${}_3F_2(1)$ summations}
\label{sec:3F2}

If $\angliceRe(c-a-b)>0$ so that ${}_2F_1(a,b;c;1)$ is convergent, its sum
can be evaluated as the gamma quotient
$\Gamma(c)\Gamma(c-a-b)/\Gamma(c-a)\Gamma(c-b)$ by Gauss's theorem.  There
is no simple analogue of this fact for ${}_{r+1}F_r(1)$, $r>1$.
Characterizing the hypergeometric series that are summable in finite terms
is an unsolved problem.

In this section the relation between the $r=2$ case of Slater's formula for
${}_{r+1}F_r(1)$ and several $3$-parameter summations with linear
parametric restrictions is investigated, and two additional $3$-parameter
${}_3F_2(1)$ summations with nonlinear parametric restrictions are derived
(see Theorem~\ref{thm:newguys}).  Underlying this section is the following
question about the level $r=2$, which at~present is too hard to answer.

{\em Which algebraic varieties $\mathcal{V}\subset
\mathbb{C}^3\times\mathbb{C}^2$ have the\/ {\em gamma quotient summation
property,} in the sense that for all parameter vectors
$(\alpha;\beta)\in\mathcal{V}$, the hypergeometric series
${}_3F_2(\alpha_1,\alpha_2,\alpha_3;\beta_1,\beta_2;1)$ sums to a gamma
quotient, the arguments of which are affine functions {\rm(}with rational
coefficients\/{\rm)} of $\alpha_1,\alpha_2,\alpha_3;\beta_1,\beta_2$?\/}
This is subject of~course to
$\angliceRe(\beta_1+\beta_2-\alpha_1-\alpha_2-\alpha_3)>0$, with
no~$\beta_i$ a non-positive integer.

A similar question could be posed about any level $r>1$, but some partial
results on $r=2$ are available.  Wimp~\cite{Wimp83} and
Zeilberger~\cite{Zeilberger92} proved that $\mathbb{C}^3\times\mathbb{C}^2$
itself does not have the gamma quotient summation property.  For a recent
discussion, including remarks on higher~$r$, see~\cite{Wimp98}.  In~Wimp's
proof, a potential obstruction is associated with any line of the form
$(\alpha_1,\alpha_2,\alpha_3;\beta_1,\beta_2)=(a,b,c;d,2c)+\xi(1,1,1;1,2)$,
$\xi\in\mathbb{C}$.  In~Zeilberger's proof a similar role is played by the
line
$(\alpha_1,\alpha_2,\alpha_3;\beta_1,\beta_2)=(0,0,1;1,1)+\xi(1,1,-1;0,0)$,
$\xi\in\mathbb{C}$.  No~variety~$\mathcal{V}$ with such a line as a
subvariety can have the property.

However, there are many varieties
$\mathcal{V}\subset\mathbb{C}^3\times\mathbb{C}^2$ that permit the
evaluation of~${}_3F_2(1)$ as a single gamma quotient.  If~an upper and a
lower parameter are equal, they may be `cancelled,' reducing the
${}_3F_2(1)$ to a ${}_2F_1(1)$.  This shows the existence of $3\cdot2=6$
hyperplanes with the property.  Less trivially, each of the classical
summation formulas of Dixon, Watson, and Whipple applies when
$(\alpha;\beta)$ is restricted to a certain $3$-dimensional affine subspace
(a~`$3$-plane') \cite[\S\tusp2.3]{Slater66}.  Dixon's formula for the sum
of a well-poised ${}_3F_2(1)$ is
\begin{sizeequation}
{\small}
\label{eq:dixon}
{}_3F_2
\biggl(
\begin{array}{c}
{a,\,b,\,c} \\
{1+a-b,\,1+a-c}\\
\end{array}
\biggr)
=\Gamma
\biggl[
\begin{array}{c}
{1+\frac12a,\,1+\frac12a-b-c,\,1+a-b,\,1+a-c} \\
{1+a,\,1+a-b-c,\,1+\frac12a-b,\,1+\frac12a-c}\\
\end{array}
\biggr],
\end{sizeequation}
and Watson's and Whipple's formulas are respectively
\begin{sizeequation}
{\small}
\label{eq:watson}
{}_3F_2\biggl(
\begin{array}{c}
{a,\,b,\,c} \\
{\frac12+\frac12a+\frac12b,\,2c}\\
\end{array}
\biggr)=\Gamma
\biggl[
\begin{array}{c}
{\frac12,\,\frac12+c,\,\frac12+\frac12a+\frac12b,\frac12-\frac12a-\frac12b+c}
\\
{\frac12+\frac12a,\,\frac12+\frac12b,\,\frac12-\frac12a+c,\frac12-\frac12b+c}
\\
\end{array}
\biggr],
\end{sizeequation}
\begin{sizeeqnarray}
{\small}
\label{eq:whipple}
&&{}_3F_2\biggl(
\begin{array}{c}
{a,\,1-a,\,c} \\
{e,\,1+2c-e}\\
\end{array}
\biggr)
=\Gamma\biggl[
\begin{array}{c}
{
e,\,1+2c-e,\,\frac12+\frac12a+\frac12e,\,\frac12-\frac12a+c-\frac12e} \\
{a+e,\,
\frac12+\frac12a+c-\frac12e,\,\frac12-\frac12a+\frac12e,
\,1-a+2c-e}\\
\end{array}
\biggr].
\end{sizeeqnarray}
The parametric excess of each~${}_3F_2(1)$ is assumed to have positive real
part, with no lower parameter equal to a non-positive integer.  The latter
two formulas are often given in alternative forms, which follow from the
duplication formula
$\Gamma(2z)=(2\pi)^{-1/2}2^{2z-1/2}\Gamma(z)\Gamma(z+1/2)$.

By Theorem~\ref{thm:slater}, the variety comprising all
$(\alpha_1,\alpha_2,\alpha_3;{\beta_1+1},{\beta_2+1})$ where
$(\alpha;\beta)\in\mathcal{U}_2$ also has the gamma quotient summation
property.  In~the absence of a list of subvarieties of
$\mathbb{C}^3\times\mathbb{C}^2$ with this property, it~is natural to
inquire as to the relation between this $3$-dimensional `Slater variety,'
which is not a $3$-plane, and the better-known $3$-planes of Dixon, Watson,
and Whipple.  Does~it play some role, perhaps, in bridging between them?

What it means for two varieties to have a high degree of overlap is
clarified by the following standard facts from algebraic geometry.  An
algebraic variety in the parameter space
$\mathbb{C}^{r+1}\times\mathbb{C}^r$ is the common zero-set of a collection
of polynomials in~$(\alpha;\beta)$.  Any irreducible variety (one with a
single component) has a well-defined dimension in the sense of algebraic
geometry; and the complement of its set of singular points, if~any, is a
Riemann surface of that dimension.  If $\mathcal{X},\mathcal{Y}$ are
varieties with dimensions $m,n$ respectively, with $m+n\ge 2r+1$, then each
component of the intersection $\mathcal{X}\cap \mathcal{Y}$, if it is
non-empty, will have dimension greater than or equal to $m+n-(2r+1)$.  If
``greater than or equal~to'' cannot be replaced by ``equal~to'', then
$\mathcal{X},\mathcal{Y}$ are said to {\em intersect improperly\/}.  Proper
intersections are generic; improper ones are~not, and are more delicate.

Suppose $\mathcal{X},\mathcal{Y}$ are $3$-dimensional varieties
in~$\mathbb{C}^3\times\mathbb{C}^2$ with the gamma quotient summation
property.  If $\mathcal{X}\cap\mathcal{Y}$ is nonempty, it will have
dimension unity and be an algebraic curve (if~the intersection is proper),
or contain a component of dimension greater than unity (if~it is improper).
The latter would signal the existence of a connection, not necessarily
obvious, between the corresponding summation formulas, since it would
indicate an unexpected amount of overlap.

It would seem easy to examine the parametric constraints of Slater's
formula and Dixon's formula, say, to determine whether when taken together,
they result in there being one free parameter or two.  But a full
examination of the extent of overlap between two summation formulas must
take into account the symmetry group of ${}_3F_2(1)$ evaluations, of which
the following is a description.

There is a trivial symmetry group isomorphic to $S_3\times S_2$, which
arises from separate permutations of upper and lower parameters.  So,
associated to any summation formula are up~to $3!\,2!=12$ trivially
equivalent formulas.  It~has long been known that there is a larger
symmetry group $\mathfrak{T}\subset GL(5,\mathbb{C})$, generated by the
permutations and by certain transformations developed by
Thomae~\cite{Thomae1879}, which linearly `mix~up' the two sorts of
parameter.  As Thomae and others showed~\cite{Hardy23}, the extended
group~$\mathfrak{T}$, which may be called the Thomae group, is isomorphic
to~$S_5$.  It~can be realized very concretely~\cite{Beyer87}.  Let
\begin{equation}
{}_3E_2
\biggl(
\begin{array}{c}
{x,\,y,\,z} \\
{u,\,v}\\
\end{array}
\biggr)
=
\frac1{\Gamma(d)\,\Gamma(e)\,\Gamma(s)}
\ {}_3F_2
\biggl(
\begin{array}{c}
{a,\,b,\,c} \\
{d,\,e}\\
\end{array}
\biggr),
\end{equation}
where $(a,b,c;d,e)=(x,y,z;u,v)A$, with the $5\times5$ matrix~$A$ given by
\begin{equation}
\label{eq:lineartrans}
A=
\left(
\begin{array}{ccccc}
1 & 0 & 0 & 1 & 1 \\
0 & 1 & 0 & 1 & 1 \\
0 & 0 & 1 & 1 & 1 \\
1 & 1 & 1 & 2 & 1 \\
1 & 1 & 1 & 1 & 2 \\
\end{array}
\right),
\qquad
3A^{-1}=
\left(
\begin{array}{rrrrr}
1 & -2 & -2 & 1 & 1 \\
-2 & 1 & -2 & 1 & 1 \\
-2 & -2 & 1 & 1 & 1 \\
1 & 1 & 1 & 1 & -2 \\
1 & 1 & 1 & -2 & 1 \\
\end{array}
\right),
\end{equation}
so that separate permutations of $x,y,z$ and~$u,v$ correspond to those of
$a,b,c$ and~$d,e$.  This defines ${}_3E_2(x,y,z;u,v)$ at all $(x,y,z;u,v)$
for which the corresponding ${}_3F_2(a,b,c;d,e;1)$ is convergent.  The
statement of Thomae-covariance is that ${}_3E_2$ extends by analytic
continuation to a function that is invariant under all $5!$ permutations
of its parameters.

Since there are $5!\,/\,3!\,2!=10$ cosets of $S_3\times S_2$ in~$S_5$,
${}_3F_2(a,b,c;d,e;1)$ can be written in $10$ distinct ways, each of the
form $G\,{}_3F_2(a',b',c';d',e';1)$, where $G$ is a gamma quotient, and the
transformed parameters $(a',b',c';d',e')$ and the arguments of~$G$ depend
linearly and homogeneously on~$(a,b,c;d,e)$.  One of these expressions is
${}_3F_2(a,b,c;d,e;1)$ itself.  The others can be extracted, with effort,
from tables originally prepared by Whipple~\cite[\S\tusp4.3]{Slater66}.
Prudnikov et~al.\ list only two~\cite[\S\tusp7.4.4, Eqs.\
(1)--(2)]{Prudnikov86c}.

If ${}_3F_2(a,b,d;d,e;1)$ is defined, i.e., if $\angliceRe(d+e-a-b-c)>0$
and neither of $d,e$ is a non-positive integer, it~is not necessarily the
case that each of the $10$~distinct expressions for~it is also defined,
since $(a',b',c';d',e')$ may fail to satisfy these conditions.  Conversely,
${}_3F_2(a',b',c';d',e';1)$ may be defined even if ${}_3F_2(a,b,c;d,e;1)$
is divergent, in which case $G\,{}_3F_2(a',b',c';d',e';1)$ may optionally
be viewed as its formal sum: which is meaningful in the sense of analytic
continuation.

\begin{theorem}
If a variety\/ $\mathcal{V}\subset\mathbb{C}^3\times\mathbb{C}^2$ has the
gamma quotient summation property, then for all\/ $5!$~transformations
$T\in\mathfrak{T}\subset GL(5,\mathbb{C})$, so does the linearly
transformed variety~$T\mathcal{V}$.  That~is, the property belongs to an
entire\/ $\mathfrak{T}$-orbit.
\end{theorem}

\begin{proof}
This is immediate.  A gamma quotient is associated to each~$T$, but a
product of gamma quotients is a gamma quotient.
\end{proof}

When investigating the relation between Slater's formula and those of
Dixon, Watson, and Whipple, one should look not merely at the intersection
of the corresponding varieties, but also at the intersections of their
images under all $T\in\mathfrak{T}$.  An~unexpected overlap could show~up
among the latter, even if the varieties themselves intersect properly.

It is well known that the Dixon, Watson, and Whipple $3$-planes lie on a
single $\mathfrak{T}$-orbit \cite[Theorem 3.5.5]{Andrews99}.  In~fact
Whipple derived his formula (and the general case of Watson's formula) by
applying appropriate $T\in\mathfrak{T}$ to Dixon's
formula~\cite{Whipple25}.  An obvious question is whether there are other
$3$-planes on the orbit.  There could be as many as~$5!$ in~all, but the
upper bound shrinks to~$10$ if $3$-planes related by permutations in
$S_3\times S_2$ are identified.  Whipple asserted without proof what
amounts to the following.

\begin{theorem}
\label{thm:DWWorbit}
Up~to separate permutations of upper and lower parameters, the
Dixon--Watson--Whipple\/ $\mathfrak{T}$-orbit in the space of\/ $3$-planes
in\/ $\mathbb{C}^3\times\mathbb{C}^2$ comprises the Dixon, Watson, and
Whipple\/ $3$-planes, and no others.
\end{theorem}

\begin{remark}
A recent paper~\cite{Exton99} produces a $3$-plane different from any of
these by applying a transformation in~$\mathfrak{T}$ to Dixon's formula,
but it unfortunately contains an algebra error.  Correcting the error
reveals that up~to separate permutations, the $3$-plane is unchanged.
\end{remark}

\begin{proof}[Proof of Theorem~\ref{thm:DWWorbit}]
The Dixon $3$-plane is specified by $d=a-b+1$ and $e=a-c+1$, or
equivalently, $(2y+z)+(2u+v)=1$ and $(y+2z)+(u+2v)=1$.  Due to the
symmetrical role played by $y,z$ and~$u,v$, its possible images under
permutations of $x,y,z,u,v$ are bijective with the partitions of
$x,y,z,u,v$ of type $2+2+1$.  Up~to separate permutations of $x,y,z$
and~$u,v$, there are only three such partitions, namely
$\{y,z\},\{u,v\},\{x\}$; and $\{y,z\},\{x,u\},\{v\}$; and
$\{y,v\},\{x,u\},\{z\}$.  By~examination, these correspond to the
$3$-planes of (\ref{eq:dixon})--(\ref{eq:whipple}).
\end{proof}

Other $3$-planes in~$\mathbb{C}^3\times\mathbb{C}^2$ with the gamma
quotient summation property are known, and when investigating the overlap
between the $r=2$ case of Slater's formula and known results, one should
consider them as~well.  No~careful description of their behavior under
Thomae's transformations seems to have been published.  This gap is
partially filled by the following theorem, which is proved similarly.

\begin{theorem}
\label{thm:PS}
For each integer $n\ge1$, there is a\/ $\mathfrak{T}$-orbit of\/ $3$-planes
in\/ $\mathbb{C}^3\times\mathbb{C}^2$ with the gamma quotient summation
property, which up~to separate permutations of upper and lower parameters
comprises exactly five\/ $3$-planes.  Two correspond to
\begin{enumerate}
\item
\label{enum:ps0}
\begin{displaymath}
{}_3F_2
\biggl(
\begin{array}{c}
{-n,\,a,\,b} \\
{c,\,1+a+b-c-n}\\
\end{array}
\biggr)
=\Gamma
\biggl[
\begin{array}{c}
{c-a+n,\,c-b+n,\,c,\,c-a-b} \\
{c-a,\,c-b,\,c+n,\,c-a-b+n} \\
\end{array}
\biggr],
\end{displaymath}
\item
\label{enum:ps1}
\begin{displaymath}
{}_3F_2
\biggl(
\begin{array}{c}
{a,\,b,\,c} \\
{a-n,\,1+b}\\
\end{array}
\biggr)
=\Gamma
\biggl[
\begin{array}{c}
{1+b-a+n,\,1-a,\,1-c,\,1+b} \\
{1+b-a,\,1-a+n,\,1+b-c}\\
\end{array}
\biggr].
\end{displaymath}
\end{enumerate}
In\/ \ref{enum:ps0} and~\ref{enum:ps1} it is assumed that no lower
parameter is a non-positive integer, and in~\ref{enum:ps1}{\rm,} that the
parametric excess has positive real part.  The remaining $3$-planes, three
in~number\/ {\rm(}up~to separate permutations\/{\rm)}, consist of points
that violate these constraints, so the corresponding identities are formal
sums of divergent ${}_3F_2(1)$'s.
\end{theorem}

\begin{remark}
The identity~\ref{enum:ps0} is the Pfaff--Saalsch\"utz formula, which sums
any terminating $1$-balanced ${}_3F_2(1)$.  It and~\ref{enum:ps1} appear in
Prudnikov et~al.\ as Eqs.\ (88),(15).  To~obtain the quotient on the
right-hand side of~\ref{enum:ps1} in the above form, the reflection formula
$\Gamma(z)\Gamma(1-z)=\pi\csc(\pi z)$ must be used repeatedly.

The remaining three identities, not given here, can be made sense~of by
analytic continuation in~$\mathbb{C}^3\times\mathbb{C}^2$
(cf.~\cite{Wimp81,Wimp83}).  The term `Pfaff--Saalsch\"utz orbit of
index~$n$' will refer to the $\mathfrak{T}$-orbit of the theorem, including
the `divergent' $3$-planes.
\end{remark}

The Pfaff--Saalsch\"utz orbits are not the only $\mathfrak{T}$-orbits of
$3$-planes, other than the Dixon--Watson--Whipple orbit, that have the
gamma quotient summation property.  An~additional orbit arises from the
summation formula
\begin{equation}
\label{eq:shukla}
{}_3F_2
\biggl(
\begin{array}{c}
{a,\,1+\frac12a,\,c} \\
{\frac12a,\,e}\\
\end{array}
\biggr)
=\Gamma
\biggl[
\begin{array}{c}
{e,\,e-a-c-1,\,e+c-a} \\
{e-a,\,e-c,\,e+c-a-1}\\
\end{array}
\biggr].
\end{equation}
This little-known identity is due to Shukla~\cite{Shukla58}, who derived it
in an indirect way.  The special terminating case (i.e., $c=-n$) was
earlier derived by Bailey, and is the only case listed in Prudnikov et~al.\
\cite[\S\tusp7.4.4, Eq.~(106)]{Prudnikov86c}.  Actually,
(\ref{eq:shukla})~has a straightforward proof.  The ${}_3F_2(1)$ has an
upper parameter that exceeds a lower by unity, so it can be written as a
combination of two ${}_2F_1(1)$'s~\cite[\S\tusp 5.2.4]{Luke75}, and
evaluated.  The $3$-plane of~(\ref{eq:shukla}) lies on its own
$\mathfrak{T}$-orbit.

Yet another $3$-parameter summation of ${}_3F_2(1)$ in which the right side
is a gamma quotient appears in Prudnikov et~al.\
\cite[\S\tusp7.4.4]{Prudnikov86c}, as Eq.~(19). Unfortunately, study
reveals that this formula, of unclear origin, is erroneous: it~is valid
only if the series terminates.  So~it does not give rise to a further
$\mathfrak{T}$-orbit of $3$-planes.

Finally, the distinct Thomae-transformed Slater summation formulas can be
enumerated.  There are only three, due~to the high degree of symmetry of
Theorem~\ref{thm:slater}.

\begin{theorem}
\label{thm:newguys}
There is a\/ $\mathfrak{T}$-orbit of\/ $3$-dimensional algebraic varieties
in\/ ${\mathbb{C}^3\times\mathbb{C}^2}$ with the gamma quotient summation
property, which up~to separate permutations of upper and lower parameters
comprises exactly three varieties, which are not\/ $3$-planes.  They
correspond to the following quadratically constrained identities.
\begin{enumerate}
\item If $ab+bc+ca = (d-1)(e-1)$ and $d+e-a-b-c=2$, then
\label{enum:slater0}
\begin{displaymath}
{}_3F_2
\biggl(
\begin{array}{c}
{a,\,b,\,c} \\
{d,\,e}\\
\end{array}
\biggr)
=\Gamma
\biggl[
\begin{array}{c}
{d,\,e} \\
{a+1,\,b+1,\,c+1}\\
\end{array}
\biggr].
\end{displaymath}
\item If $(a-1)(b-1)=[(a-1)+(b-1)-(e-1)]c$, then
\label{enum:slater1}
\begin{displaymath}
{}_3F_2
\biggl(
\begin{array}{c}
{a,\,b,\,c} \\
{c+2,\,e}\\
\end{array}
\biggr)
=\Gamma
\biggl[
\begin{array}{c}
{e,\,e-a-b+2,\,c+2} \\
{e-a+1,\,e-b+1,\,c+1}\\
\end{array}
\biggr],
\end{displaymath}
provided\/ $\angliceRe(e-a-b+2)>0$.
\item If $(a-1)(b-1)=(d-2)(e-2)$, then
\label{enum:slater2}
\begin{displaymath}
{}_3F_2
\biggl(
\begin{array}{c}
{a,\,b,\,2} \\
{d,\,e}\\
\end{array}
\biggr)
=\frac{(d-1)(e-1)}{d+e-a-b-2},
\end{displaymath}
provided\/ $\angliceRe(d+e-a-b-2)>0$.
\end{enumerate}
In each identity, it is assumed that no lower parameter is a
non-positive integer.
\end{theorem}

\begin{remark}
Each of \ref{enum:slater0}--\ref{enum:slater2}, if the quadratic constraint
on its parameters is dropped, can be generalized to a formula in which a
second ${}_3F_2(1)$ appears on the right-hand side.  For example,
\ref{enum:slater0}~generalizes to a non-terminating version of the
Sheppard--Anderson identity, and \ref{enum:slater2}~to Eq.~(27) of
Prudnikov et~al.\ \cite[\S\tusp7.4.4]{Prudnikov86c}.

Numerical experimentation reveals that \ref{enum:slater1}
and~\ref{enum:slater2} hold whenever the series terminates, even if its
parametric excess has a non-positive real part.  In~this regard they differ
from the formulas of Dixon, Watson, and Whipple.  In~the general
terminating case, those three formulas must be
modified~\cite{Bailey53,Dzhrbashyan64}.  

If a lower parameter is a non-positive integer and an upper one is a
non-positive integer of smaller magnitude, the series is usually regarded
as terminating.  Experimentation reveals that \ref{enum:slater2}~holds in
this case, though \ref{enum:slater0} and~\ref{enum:slater1} do not.
\end{remark}

\begin{proof}[Proof of Theorem~\ref{thm:newguys}]
Formula~\ref{enum:slater0} is a restatement of the $r=2$ case of
Theorem~\ref{thm:slater}.  It~follows from~(\ref{eq:lineartrans}) that
in~terms of $x,y,z,u,v$, its variety is
\begin{eqnarray*}
\left\{
\begin{array}{rcl}
xy+yz+zx + u^2 + v^2 + uv + u + v &=& 1 \\
x+y+z&=&2.\\
\end{array}
\right.
\end{eqnarray*}
Due to the symmetrical role played by $x,y,z$ and~$u,v$, up~to separate
permutations there are only three images of this variety under a
permutation $\sigma\in S_5$.  It~suffices to consider the permutations
$x,y,z,u,v$; and $x,y,u,z,v$; and $x,v,u,z,y$.  By~examination, these yield
\ref{enum:slater0}--\ref{enum:slater2}, respectively.  Interestingly, the
gamma quotient in~\ref{enum:slater2} reduces to a rational function of the
parameters.
\end{proof}

Just as the $r=2$ case of Theorem~\ref{thm:slater} gives rise to the
$\mathfrak{T}$-orbit of Theorem~\ref{thm:newguys}, the $r=2$ case of
Theorem~\ref{thm:slaterplus} gives rise to a $\mathfrak{T}$-orbit of
$3$-dimensional algebraic varieties.  However, the three ${}_3F_2(1)$
summations on that orbit are not of gamma quotient type.  It~is clear from
the right side of Theorem~\ref{thm:slaterplus} that each yields not a gamma
quotient, but a gamma quotient times a rational function of the
hypergeometric parameters.

The extent to which the identities of Theorem~\ref{thm:newguys} and those
of Dixon, Watson, Whipple and Pfaff--Saalsch\"utz overlap one another will
now be considered.  The term `Slater orbit' will refer to the 
$\mathfrak{T}$-orbit of Theorem~\ref{thm:newguys}.

\begin{definition}
If $\mathfrak{V}$ and $\mathfrak{W}$ are two $\mathfrak{T}$-orbits of
algebraic varieties in $\mathbb{C}^3\times\mathbb{C}^2$ (satisfying
$m+n\ge5$, where $m,n$ are the dimensions of the varieties in
$\mathfrak{V},\mathfrak{W}$ respectively), they are said to intersect
improperly iff there is at~least one $\mathcal{V}\in\mathfrak{V}$ and one
$\mathcal{W}\in\mathfrak{W}$ that intersect improperly.  By~invariance
under~$\mathfrak{T}$, this occurs iff for all $\mathcal{V}\in\mathfrak{V}$,
there is at~least one $\mathcal{W}\in\mathfrak{W}$ such that $\mathcal{V}$
and~$\mathcal{W}$ intersect improperly.
\end{definition}

\begin{theorem}
\label{thm:orbits}
\ 
\begin{enumerate}
\item The Slater orbit and Dixon--Watson--Whipple orbit intersect improperly.
\label{enum:slaterdww}
\item 
\label{enum:slaterps}
For all\/ $n\ge1$, the Slater orbit and the Pfaff--Saalsch\"utz orbit of
index~$n$ do not intersect improperly.
\item 
\label{enum:dwwps}
For all\/ $n\ge1$, the Dixon--Watson--Whipple orbit and the
Pfaff--Saal\-sch\"utz orbit of index~$n$ do not intersect improperly.
\end{enumerate}
\end{theorem}

\begin{proof}
\ref{enum:slaterdww}.  It suffices to show that the parametric restrictions
of Theorem~\ref{thm:newguys}\ref{enum:slater0} and Dixon's formula, taken
together, result in there being two free parameters rather than one.  By
examination, the $2$-parameter formula
\begin{equation}
\label{eq:overlap1}
{}_3F_2
\biggl(
\begin{array}{c}
{2b+2c,\,b,\,c} \\
{1+b+2c,\,1+2b+c}\\
\end{array}
\biggr)
=\Gamma
\biggl[
\begin{array}{c}
{1+b+2c,\,1+2b+c} \\
{1+2b+2c,\,1+b,\,1+c}\\
\end{array}
\biggr],
\end{equation}
the variety corresponding to which is a $2$-plane, is a specialization of
both Dixon's formula and Theorem~\ref{thm:newguys}\ref{enum:slater0}; so
their overlap is improper.

\ref{enum:slaterps}.  It suffices to compare
Theorem~\ref{thm:PS}\ref{enum:ps0} or Theorem~\ref{thm:PS}\ref{enum:ps1}
with each of
Theorem~\ref{thm:newguys}\ref{enum:slater0},\ref{enum:slater1},\ref{enum:slater2},
and verify that their parametric constraints, taken together, result in
there being a single free parameter.  (Each comparison is really $3!\,2!$
separate ones, since all separate permutations of upper and lower
parameters must be considered.)  It~is easier to use
Theorem~\ref{thm:PS}\ref{enum:ps0}.  Its series is $1$-balanced and that of
Theorem~\ref{thm:newguys}\ref{enum:slater0} is $2$-balanced, so the
corresponding varieties do not intersect; and only
Theorem~\ref{thm:newguys}\ref{enum:slater1},\ref{enum:slater2} need to be
compared with.  Details are left to the reader.

\ref{enum:dwwps}.  This is similar to the proof of~\ref{enum:slaterps}.
It~is easy to compare Theorem~\ref{thm:PS}\ref{enum:ps1} with each of
(\ref{eq:dixon})--(\ref{eq:whipple}), and verify that each possible overlap
(separate permutations of upper and lower parameters being allowed) has
only one free parameter.
\end{proof}

It follows from Theorem~\ref{thm:orbits}\ref{enum:slaterdww} that each of
the Dixon, Watson, and Whipple $3$-planes improperly intersects at~least
one of the three algebraic varieties of Theorem~\ref{thm:newguys}.
Conversely, each of those algebraic varieties improperly intersects
at~least one of the Dixon, Watson, and Whipple $3$-planes.  In~these
statements the phrase ``up~to separate permutations of upper and lower
parameters'' is understood.

The following examples may be instructive.  Up~to separate permutations,
\begin{equation}
\label{eq:overlap2}
{}_3F_2
\biggl(
\begin{array}{c}
{2,\,b,\,c} \\
{\frac32+\frac12b,\,2c}\\
\end{array}
\biggr) = \frac {(1+b)(1-2c)} {(1+b-2c)}
\end{equation}
is a specialization of both Watson's formula and
Theorem~\ref{thm:newguys}\ref{enum:slater2}.  Similarly,
\begin{equation}
\label{eq:overlap3}
{}_3F_2
\biggl(
\begin{array}{c}
{a,\,1-a,\,c} \\
{2+a,\,2c-a-1}\\
\end{array}
\biggr)
=\Gamma
\biggl[
\begin{array}{c}
{2c-a-1,\,c,\,2+a} \\
{c-a,\,2c-1,\,1+a}\\
\end{array}
\biggr]
\end{equation}
is a specialization of both Whipple's formula and
Theorem~\ref{thm:newguys}\ref{enum:slater1}.  Additional examples of
$2$-parameter overlap can be constructed.
In~$\mathbb{C}^3\times\mathbb{C}^2$, the Slater variety and its images
under Thomae's transformations intersect the Dixon--Watson--Whipple
$3$-planes quite extensively.

The preceding results, such as (\ref{eq:overlap1})--(\ref{eq:overlap3}),
clarify a remark of Slater \cite[\S\tusp2.6.1]{Slater66}
that the $r=2$ case of her summation formula is ``a~disguised form of
Dixon's theorem''.  It~is better described as an extension of a
$2$-parameter special case of Dixon's theorem; and if Thomae's
transformations are taken into account, the unexpected degree of overlap
between the two summation formulas becomes an unexpected degree of
overlap between two $\mathfrak{T}$-orbits.

\section{${}_2F_1(-1)$ summations}
\label{sec:2F1}

If $\angliceRe(c-a-b)>-1$ then ${}_2F_1(a,b;c;-1)$ is convergent, though in
general it cannot be evaluated in~terms of the gamma function.  In~the
well-poised case, it can be evaluated as a single gamma quotient.  Kummer's
theorem states that
\begin{equation}
\label{eq:kummer}
{}_2F_1
\biggl(
\begin{array}{c}
{a,\,b} \\
{1+a-b}\\
\end{array}
\!
\biggm|
-1
\biggr)
=\Gamma
\biggl[
\begin{array}{c}
{1+a-b,\,1+\frac12a} \\
{1+a,\,1+\frac12a-b}\\
\end{array}
\biggr],
\end{equation}
provided $1+a-b$ and $1+\frac12a$ are not non-positive integers.  Prudnikov
et~al.\ \cite[\S\tusp7.3.6]{Prudnikov86c} supply a list of known
evaluations of~${}_2F_1(-1)$ with one or more free parameters.  Besides
Kummer's formula, these include several multi-term evaluations contiguous
to~it, in the sense that they follow by contiguous function relations.


Kummer's formula is often viewed as a consequence of a quadratic
transformation on the ${}_2F_1$ level, based on
$x\mapsto
R(x)=-4x/(1-x)^2$, which applies to well-poised
series~\cite[\S\tusp3.1]{Andrews99}.  Since $R(-1)=1$, the transformation
permits the summation of well-poised ${}_2F_1(-1)$ series by Gauss's
formula, leading to~(\ref{eq:kummer}).  Kummer's formula will be viewed
differently here.
Whipple~\cite{Whipple29} showed that
\begin{equation}
\label{eq:whipple2f1}
{}_2F_1
\biggl(
\begin{array}{c}
{2\alpha,\,\beta} \\
{2\gamma-\beta}\\
\end{array}
\!
\biggm|
-1
\biggr)
=\Gamma
\biggl[
\begin{array}{c}
{2\gamma-2\alpha,\,2\gamma-\beta} \\
{2\gamma-2\alpha-\beta,\,2\gamma}\\
\end{array}
\biggr]
{}_3F_2
\biggl(
\begin{array}{c}
{\alpha,\,\frac12+\alpha,\,\beta} \\
{\gamma,\,\frac12+\gamma}\\
\end{array}
\biggr),
\end{equation}
provided that both sides are defined.  So there is a $3$-plane in the
$\mathbb{C}^3\times\mathbb{C}^2$ parameter space, i.e.,
$\mathcal{W}_0=\{\,(\alpha_1,\alpha_2,\alpha_3;\beta_1,\beta_2)\mid
\alpha_2=\alpha_1+\frac12,\,\beta_2=\beta_1+\frac12\,\}$, on which any
evaluation of ${}_3F_2(1)$ leads to an evaluation of a~${}_2F_1(-1)$.  By
covariance under the Thomae group~$\mathfrak{T}$, the same is true of all
$3$-planes of the form $T\mathcal{W}_0$, $T\in\mathfrak{T}$.  Whipple
showed that up~to separate permutations of upper and lower parameters,
there are exactly six such $3$-planes, including $\mathcal{W}_0$ itself,
and the term `Whipple ${}_2F_1(-1)$ orbit' will refer to this
$\mathfrak{T}$-orbit.  His paper contains five additional
Thomae-transformed versions of the relation~(\ref{eq:whipple2f1}), which
will not be needed explicitly.

If any $\mathfrak{T}$-orbit of algebraic varieties
in~$\mathbb{C}^3\times\mathbb{C}^2$ has the gamma quotient summation
property, its intersection with the Whipple ${}_2F_1(-1)$ orbit will too,
if non-empty, and one or more evaluations of ${}_2F_1(-1)$ as a single
gamma quotient will result.  In~the framework of this paper, this is the
origin of~(\ref{eq:kummer}).  {\em Kummer's formula arises from the
non-empty intersection of the Dixon--Watson--Whipple orbit with the Whipple
${}_2F_1(-1)$ orbit.}  In~particular, it arises from the intersection of a
Dixon $3$-plane with the Whipple ${}_2F_1(-1)$ $3$-plane, which produces a
$2$-plane rather than a line (and~is therefore improper).  Combining
(\ref{eq:dixon}) with~(\ref{eq:whipple2f1}) yields the $2$-parameter
formula
\begin{align*}
&{}_2F_1
\biggl(
\begin{array}{c}
{a,\,b} \\
{1+a-b}\\
\end{array}
\!
\biggm|
-1
\biggr)
\\
&\qquad=
\Gamma
\biggl[
\begin{array}{c}
{1+2a-2b,\,1+a-b} \\
{1+a-2b,\,1+2a-b}\\
\end{array}
\biggr]
\,{}_3F_2
\biggl(
\begin{array}{c}
{\frac12b,\,\frac12+\frac12b,\,a} 
\\
{\frac12+a-\frac12b,\,1+a-\frac12b}\\
\end{array}
\biggr) 
\\
&\qquad=
\Gamma
\biggl[
\begin{array}{c}
{1+2a-2b,\,1+a-b} \\
{1+a-2b,\,1+2a-b}\\
\end{array}
\biggr]
\\
&\qquad\qquad\qquad\qquad
{}\times
\Gamma
\biggl[
\begin{array}{c}
{1+\frac12a,\,\frac12+\frac12a-b,\,1+a-\frac12b,\,\frac12+a-\frac12b} \\
{1+a,\,\frac12+a-b,\,1+\frac12a-\frac12b,\,\frac12+\frac12a-\frac12b} \\
\end{array}
\biggr]
\end{align*}
If this is simplified with the aid of the duplication formula for the gamma
function, it becomes Kummer's formula.

In the same way, many ${}_2F_1(-1)$ summations can be derived by
intersecting the Pfaff--Saalsch\"utz and Slater orbits with the Whipple
${}_2F_1(-1)$ orbit.  Each intersection turns~out to be proper, producing a
union of lines, or (in~the Slater case) algebraic curves.  The
$1$-parameter summation formulas corresponding to the latter, with
nonlinear parametric constraints, are quite exotic.  They are given in
Theorem~\ref{thm:nonlinear} and (\ref{eq:eulered1})--(\ref{eq:eulered2})
below.

Actually, there are intersections between the Dixon--Watson--Whipple and
Whipple ${}_2F_1(-1)$ orbits other than the improper one yielding Kummer's
formula.  Surprisingly, the resulting $1$-parameter identities have never
been systematically worked out.  Theorem~\ref{thm:kummercontiguous} does
this.  It~refers to Euler's transformation in the form
\begin{multline}
\label{eq:neweuler}
{}_2F_1\biggl(
\begin{array}{c}
{a,\,b} \\
{c}\\
\end{array}
\!
\biggm|
-1
\biggr)
=
\Gamma
\biggl[
\begin{array}{c}
{\frac12,\,1+c-a-b} \\
{\frac12+\frac12c-\frac12a-\frac12b,\,1+\frac12c-\frac12a-\frac12b}\\
\end{array}
\biggr]
\\
{}\times
{}_2F_1\biggl(
\begin{array}{c}
{c-a,\,c-b} \\
{c}\\
\end{array}
\!
\biggm|
-1
\biggr),
\end{multline}
which follows from~(\ref{eq:euler}) by the duplication formula, and leaves
the gamma quotient summation property for ${}_2F_1(-1)$ unaffected.  Note
that in~(\ref{eq:neweuler}) the parametric excess~$s$ of each series is the
negative of that of the other.  Since each series, if non-terminating,
converges only if $\angliceRe s>-1$, it~is possible for one to be
convergent and the other divergent.  If so, (\ref{eq:neweuler})~will supply
a formal sum for the divergent one.

The reason why (\ref{eq:neweuler}) is relevant here is the following.  The
Whipple ${}_2F_1(-1)$ orbit comprises all $3$-planes of the form
$T\mathcal{W}_0$, $T\in\mathfrak{T}$.  Whipple's
identity~(\ref{eq:whipple2f1}) supplies a canonical map $P$
of~$\mathcal{W}_0$ onto $\mathbb{C}^2\times\mathbb{C}^1$, the parameter
space of ${}_2F_1(-1)$.  Any $3$-plane~$\mathcal{W}=T\mathcal{W}_0$ is
projected onto $\mathbb{C}^2\times\mathbb{C}^1$ by~$PT^{-1}$.  This map is
not unique, since there are Thomae transformations~$T$ other than the
identity that map $\mathcal{W}_0$ onto itself, yet have the property that
$PT^{-1}\neq P$.  But as Whipple noted in somewhat different
language~\cite{Whipple29}, the map from each~$\mathcal{W}$
onto~$\mathbb{C}^2\times\mathbb{C}^1$ is unique up~to composition with a
subsequent map $(a,b;c)\mapsto(c-a,c-b;c)$, i.e., up~to a subsequent
Euler's transformation.

\begin{theorem}
\label{thm:kummercontiguous}
Up to separate permutations of upper and lower parameters and up~to Euler's
transformation, the intersection of the\/ $3$-planes on the
Dixon--Watson--Whipple\/ $\mathfrak{T}$-orbit and the\/ $3$-planes on the
Whipple ${}_2F_1(-1)$ $\mathfrak{T}$-orbit, when projected to\/
$\mathbb{C}^2\times\mathbb{C}^1$ via the map in Whipple's
identity\/~{\rm(\ref{eq:whipple2f1})}, comprises\/ {\rm(a)}~the Kummer\/
$2$-plane, and\/ {\rm(b)}~exactly four lines.  The latter correspond to the
following.
\begin{enumerate}
\item 
\label{enum:kummer1}
\begin{displaymath}
{}_2F_1\biggl(
\begin{array}{c}
{a-4,\,\frac23a-1} \\
{1+\frac13a}\\
\end{array}
\!
\biggm|
-1
\biggr)
=\Gamma
\biggl[
\begin{array}{c}
{\frac12,\,\frac13a,\,1+\frac13a,\,\frac23a-\frac32,\,\frac23a-2} \\
{\frac12a-\frac32,\,\frac16a,\,\frac12+\frac16a,\,2-\frac16a,\,\frac16a,\,\frac43a-3} \\
\end{array}
\biggr],
\end{displaymath}
\item
\label{enum:kummer2}
\begin{displaymath}
{}_2F_1\biggl(
\begin{array}{c}
{a-\frac32,\,\frac14a+\frac14} \\
{\frac34a-\frac14}\\
\end{array}
\!
\biggm|
-1
\biggr)
= 
\Gamma
\biggl[
\begin{array}{c}
{\frac12,\,\frac32,\,\frac34+\frac14a,\,\frac12a,\,\frac34a-\frac14} \\
{\frac78,\,\frac98,\,a,\,\frac18+\frac14a,\,\frac38+\frac14a}\\
\end{array}
\biggr],
\end{displaymath}
\item
\label{enum:kummer3}
\begin{displaymath}
{}_2F_1\biggl(
\begin{array}{c}
{a-\frac12,\,\frac14a} \\
{\frac34a}\\
\end{array}
\!
\biggm|
-1
\biggr)
=
\Gamma
\biggl[
\begin{array}{c}
{\frac12,\,\frac12,\,\frac12+\frac14a,\,\frac12+\frac12a,\,\frac34a} \\
{\frac38,\,\frac58,\,a,\,\frac38+\frac14a,\,\frac58+\frac14a}\\
\end{array}
\biggr],
\end{displaymath}
\item
\label{enum:kummer4}
\begin{displaymath}
{}_2F_1\biggl(
\begin{array}{c}
{a+2,\,\frac23a} \\
{\frac13a}\\
\end{array}
\!
\biggm|
-1
\biggr)
=\Gamma
\biggl[
\begin{array}{c}
{\frac12,\,\frac13a,\,1+\frac13a,\,\frac32+\frac23a,\,2+\frac23a} \\
{\frac32+\frac12a,\,\frac12+\frac16a,\,1+\frac16a,\,-\frac16a,\,2+\frac43a} \\
\end{array}
\biggr].
\end{displaymath}
\end{enumerate}
In each of the above it is assumed that the parametric excess of the
${}_2F_1(-1)$ series has real part greater than\/~$-1$, with no lower
parameter of the series or upper parameter of the gamma quotient equal to a
non-positive integer.
\end{theorem}

\begin{remark}
Each of \ref{enum:kummer1}--\ref{enum:kummer4} is written as a
${}_2F_1(a+\nu,b;a-b;-1)$ evaluation, to agree with the convention of
Vid\=unas~\cite{Vidunas2001a}.  Respectively, $\nu$~equals
$-4,-\frac32,-\frac12,2$.  In~each, $b$~is some linear function of~$a$.
The quantity $\nu$ measures the nearness of the series to well-poisedness.
When $\nu=-1$, the series is well-poised and can be summed by Kummer's
formula, with no restriction on~$b$, to give the gamma quotient
$\Gamma(a-b)\Gamma(\frac12+\frac12a)/\Gamma(a)\Gamma(\frac12+\frac12a-b)$.

Note that by using the duplication formula, the right side
of~\ref{enum:kummer1} can be written in a simpler way, as
$(3/4)\Gamma(1+\frac13a)\Gamma(\frac12a-1)/\Gamma(a-2)\Gamma(2-\frac16a)$.
The right sides of \ref{enum:kummer2}--\ref{enum:kummer4} can also be
rewritten in various ways.
\end{remark}

\begin{proof}[Proof of Theorem~\ref{thm:kummercontiguous}]
This is straightforward though tedious.  It suffices to compare the Whipple
${}_2F_1(-1)$ formula~(\ref{eq:whipple2f1}) with each of
(\ref{eq:dixon})--(\ref{eq:whipple}), and determine the possible ways in
which their respective parametric constraints on~${}_3F_2(1)$ can be
simultaneously satisfied.  (Each of these three comparisons is really
$3!\,2!$ separate comparisons.)  Every success yields a non-empty
intersection of the corresponding $3$-planes, and a parametrized gamma
quotient evaluation of~${}_2F_1(-1)$.

An~example is the following.  Specializing parameters
in~(\ref{eq:whipple2f1}) yields
\begin{displaymath}
{}_2F_1\biggl(
\begin{array}{c}
{2\alpha,\,3\alpha-\frac52} \\
{\frac32+\alpha}\\
\end{array}
\!
\biggm|
-1
\biggr)
=\Gamma
\biggl[
\begin{array}{c}
{2\alpha-1,\,\frac32+\alpha} \\
{\frac32-\alpha,\,4\alpha-1} \\
\end{array}
\biggr]
\,{}_3F_2\biggl(
\begin{array}{c}
{\alpha,\,\frac12+\alpha,\,3\alpha-\frac52} \\
{2\alpha-\frac12,\,2\alpha} \\
\end{array}
\biggr),
\end{displaymath}
in which the ${}_3F_2(1)$ may be evaluated by~(\ref{eq:watson}), Watson's
formula.  Performing the evaluation and changing the independent variable
from $\alpha$ to $a=\frac32-3\alpha$ yields the $1$-parameter evaluation of
${}_2F_1(-1)$ shown in
Theorem~\ref{thm:kummercontiguous}\ref{enum:kummer1}.

Besides the $2$-parameter formula of Kummer, the comparisons yield eight
$1$-parameter formulas in~all, but these split into four pairs, each being
related by~(\ref{eq:neweuler}), Euler's transformation.
Theorem~\ref{thm:kummercontiguous} lists only a single member of each.
\end{proof}

Theorem~\ref{thm:kummercontiguous} is closely related to the results of
Vid\=unas, who recently worked~out a systematic way of generating
${}_2F_1(-1)$ evaluations contiguous to Kummer's
formula~\cite{Vidunas2001a}.  He showed that for all integer~$\nu$,
\begin{equation}
\label{eq:vidunas}
{}_2F_1\biggl(
\begin{array}{c}
{a+\nu,\,b} \\
{a-b}\\
\end{array}
\!
\biggm|
-1
\biggr)
=
P(\nu;a,b)
\,\Gamma
\biggl[
\begin{array}{c}
{a-b,\,\frac12+\frac12a} \\
{a,\,\frac12+\frac12a-b} \\
\end{array}
\biggr]
+Q(\nu;a,b)
\,\Gamma
\biggl[
\begin{array}{c}
{a-b,\,\frac12a} \\
{a,\,\frac12a-b} \\
\end{array}
\biggr],
\end{equation}
where $P(\nu;a,b)$ and $Q(\nu;a,b)$ are rational in $a,b$.  For any
$\nu\in\mathbb{Z}$, setting $P(\nu;a,b)$ or $Q(\nu;a,b)$ to zero constrains
$a,b$ algebraically; and if the constraint is satisfied, the remaining term
on the right-hand side of~(\ref{eq:vidunas}) becomes (almost) a gamma
quotient summation of ${}_2F_1(-1)$.  If $Q(\nu;a,b)$ or $P(\nu;a,b)$,
respectively, factors into quotients of linear polynomials, then it will
be, in~fact, a gamma quotient.

Vid\=unas noticed that if $\nu=-4$ then $b$~is constrained to be a linear
function of~$a$, and he worked~out a ${}_2F_1(-1)$ summation equivalent to
Theorem~\ref{thm:kummercontiguous}\ref{enum:kummer1}.  When $\nu=2$, his
approach also yields a gamma quotient summation, equivalent to
Theorem~\ref{thm:kummercontiguous}\ref{enum:kummer4}.  A~nice feature of
his treatment is that it makes clear where the reflection symmetry through
$\nu=-1$ comes from.  (The set $\{-4,-\frac32,-\frac12,2\}$ is invariant
under this operation.)  He~showed that if $\nu\in\mathbb{Z}$, reflecting
$P(\nu;a,b)$, $Q(\nu;a,b)$ through $\nu=-1$ simply multiplies them by gamma
quotients.  However, his results do not immediately extend to
non-integer~$\nu$.

Finally, the exotic ${}_2F_1(-1)$ summations arising from Slater's formula
can be given.

\begin{theorem}
\label{thm:nonlinear}
Up to separate permutations of upper and lower parameters and up~to Euler's
transformation, the intersection of the\/ $3$-dimensional algebraic
varieties on the Slater\/ $\mathfrak{T}$-orbit and the\/ $3$-planes on the
Whipple ${}_2F_1(-1)$ $\mathfrak{T}$-orbit, when projected to\/
${\mathbb{C}^2\times\mathbb{C}^1}$ via the map in Whipple's
identity\/~{\rm(\ref{eq:whipple2f1})}, comprises\/ {\rm(a)}~three lines
that lie in the Kummer\/ $2$-plane, {\rm(b)}~one line that lies in the
displaced\/ $2$-plane with $\nu=-4$, and\/ {\rm(c)}~two algebraic curves
other than lines.  The last correspond to the following
identities, with parameters constrained both quadratically and linearly.
\begin{enumerate}
\item
If\/ $2(\alpha-\beta)^2-3\alpha+4\beta+1=0$, then
\label{enum:nonlinear1}
\begin{align*}
{}_2F_1\biggl(
\begin{array}{c}
{\beta-2\alpha+3,\,3} \\
{\beta+3}\\
\end{array}
\!
\biggm|
-1
\biggr)
=
\frac14
\left[
\frac{(\beta+1)(\beta+2)}{2\beta-2\alpha+3}
\right].
\end{align*}
\item
If\/ $2(\alpha-\beta)^2-3\alpha+6\beta+1=0$, then
\label{enum:nonlinear2}
\begin{align*}
{}_2F_1\biggl(
\begin{array}{c}
{\beta-2\alpha+4,\,4} \\
{\beta+4}\\
\end{array}
\!
\biggm|
-1
\biggr)
=\frac14
\left[
\frac{(\beta+1)(\beta+2)(\beta+3)}
{(2\beta-2\alpha+4)(2\beta-2\alpha+5)}
\right].
\end{align*}
\end{enumerate}
In both\/ \ref{enum:nonlinear1} and\/~\ref{enum:nonlinear2} it is assumed
that the ${}_2F_1(-1)$ series converges.
\end{theorem}

\begin{proof}
Like the proof of Theorem~\ref{thm:kummercontiguous}, this is
straightforward but lengthy.  It suffices to compare the Whipple
${}_2F_1(-1)$ formula~(\ref{eq:whipple2f1}) with each of
Theorem~\ref{thm:newguys}\ref{enum:slater0}--\ref{enum:slater2}, and
determine the ways in which their respective parametric constraints
on~${}_3F_2(1)$ can be simultaneously satisfied.  (Each comparison is
really $3!\,2!$ separate ones.)  Every success yields a non-empty
intersection of a Whipple ${}_2F_1(-1)$ $3$-plane with a $3$-dimensional
variety on the Slater orbit, and a parametrized gamma quotient evaluation
of~${}_2F_1(-1)$.  Each intersection turns~out to be proper, i.e.,
$1$-dimensional.

An~example is the following.  Specializing parameters
in~(\ref{eq:whipple2f1}) yields
\begin{displaymath}
{}_2F_1
\biggl(
\begin{array}{c}
{2\alpha,\,\beta} \\
{3+\beta}\\
\end{array}
\!
\biggm|
-1
\biggr)
=\Gamma
\biggl[
\begin{array}{c}
{3+2\beta-2\alpha,\,3+\beta} \\
{3+\beta-2\alpha,\,3+2\beta}\\
\end{array}
\biggr]
{}_3F_2
\biggl(
\begin{array}{c}
{\alpha,\,\frac12+\alpha,\,\beta} \\
{\frac32+\beta,\,2+\beta}\\
\end{array}
\biggr).
\end{displaymath}
If $2(\alpha-\beta)^2-3\alpha+4\beta+1=0$, the ${}_3F_2(1)$ may be
evaluated by Theorem~\ref{thm:newguys}\ref{enum:slater1}.  Performing the
evaluation, and also applying Euler's transformation in the
form~(\ref{eq:neweuler}) to the ${}_2F_1(-1)$ on the left-hand side, yields
the $1$-parameter formula for ${}_2F_1(-1)$ shown in
Theorem~\ref{thm:nonlinear}\ref{enum:nonlinear1}.
\end{proof}

By uniformizing the two algebraic curves of the theorem, i.e.,
parametrizing them by an auxiliary variable $t\in\mathbb{C}$, the
identities \ref{enum:nonlinear1} and~\ref{enum:nonlinear2} of the theorem
can be converted to the mysterious parametric forms
\begin{align}
\label{eq:eulered1}
{}_2F_1\biggl(
\begin{array}{c}
{2t^2-5t+4,\,3} \\
{-2t^2+3t+2} \\
\end{array}
\!
\biggm|
-1
\biggr)
&= -\frac14\,t\,(2t^2-3t-1),
\\[\jot]
\label{eq:eulered2}
{}_2F_1\biggl(
\begin{array}{c}
{6t^2-11t+6,\,4} \\
{-6t^2+5t+3} \\
\end{array}
\!
\biggm|
-1
\biggr)
&= -\frac1{24}\,t\,(6t+1)\,(6t^2-5t-2).
\end{align}
Each holds for all $t\in\mathbb{C}$ for which the ${}_2F_1(-1)$ series
converges.

The identities (\ref{eq:eulered1})--(\ref{eq:eulered2}) could alternatively
be derived from Gauss's three-term contiguous function relations
for~${}_2F_1$, mentioned in Section~\ref{sec:key}.  In~fact, iterating
the relation
\begin{multline*}
(c-b)\ 
{}_2F_1\biggl(
\begin{array}{c}
{a,\,b-1} \\
{c} \\
\end{array}
\!
\biggm|
x
\biggr)
+
(2b-c-bx+ax)\ 
{}_2F_1\biggl(
\begin{array}{c}
{a,\,b} \\
{c} \\
\end{array}
\!
\biggm|
x
\biggr)\\
{}+b(x-1)\ 
{}_2F_1\biggl(
\begin{array}{c}
{a,\,b+1} \\
{c} \\
\end{array}
\!
\biggm|
x
\biggr)
=0
\end{multline*}
so as to express ${}_2F_1(a,b;c;x)$ with $b=2,3$ as a combination of
${}_2F_1(a,0;c;x)\equiv1$ and ${}_2F_1(a,1;c;x)$, and setting the
coefficient of ${}_2F_1(a,1;c;x)$ to zero and uniformizing the resulting
$x$-dependent algebraic curve in $\mathbb{C}^2\ni(a,c)$, yields the
identities
\begin{equation}
{}_2F_1\biggl(
\begin{array}{c}
{-x^{-1}t + (1-x^{-1}),\,2} \\
{-t + 1} \\
\end{array}
\!
\biggm|
x
\biggr)
= \frac{t}{x-1},
\end{equation}
\begin{multline}
{}_2F_1\biggl(
\begin{array}{c}
{(1-x^{-1})t^2 + (-2+3x^{-1})t + (2-2x^{-1}),\,3} \\
{(x-1)t^2 + (-x+2)t + 2} \\
\end{array}
\!
\biggm|
x
\biggr)\\
= \frac{t\left[(-x+1)t^2 + (x-2)t -1\right]}{2(x-1)},
\end{multline}
the latter of which subsumes~(\ref{eq:eulered1}).  This sequence of
identities cannot be continued indefinitely, since the $b=4$ curve is
cubic: it~is the zero-set of a cubic polynomial in~$(a,c)$, with
coefficients polynomial in~$x$.  The genus of an irreducible algebraic
curve of degree~$n$ is $\binom{n-1}2$ minus the number of its double
points, which is generically zero.  So one would expect that generically,
i.e., except at isolated values of~$x$, this cubic curve would be
irreducible of genus~$1$, ruling~out a uniformization by rational
functions.  Some computation reveals that the situation is not so bad as
that: irrespective of~$x$, the $b=4$ curve always has genus zero.  In~fact,
if $x$~is one of $-1,1/2,2$, then it turns~out to be reducible: the union
of a quadratic curve and a line.  When $x=-1$, the quadratic component can
be uniformized by polynomials, yielding~(\ref{eq:eulered2}).  However, the
problem of positive genus does set~in when $b\ge5$.  Fuller details of this
alternative approach will appear elsewhere.

Actually, the first proof of (\ref{eq:eulered1})--(\ref{eq:eulered2}),
which was based not on contiguous function relations but rather on
Theorem~\ref{thm:nonlinear}, and therefore on the generalized Euler's
transformation Theorem~\ref{thm:main}, has much to recommend~it.  Its
existence suggests that Theorem~\ref{thm:main} encapsulates much of the
power of the contiguous function relations for ${}_{r+1}F_r$ to generate
one-term summation formulas.


\small\def\em{\it} \newcommand{\noopsort}[1]{} \newcommand{\printfirst}[2]{#1}
  \newcommand{\singleletter}[1]{#1} \newcommand{\switchargs}[2]{#2#1}
\providecommand{\bysame}{\leavevmode\hbox to3em{\hrulefill}\thinspace}

\end{document}